\documentclass[12pt]{amsart}

\usepackage{tikz}
\usetikzlibrary{positioning,automata}
\tikzset{every state/.style={minimum size=3mm}}

\usetikzlibrary{calc}

\tikzset{onslide/.code args={<#1>#2}{%
  \only<#1>{\pgfkeysalso{#2}} 
}}

\usepackage{amsmath}
\usepackage{amssymb}
\usepackage{url}
\usepackage[all,ps]{xy}
\usepackage{tikz}

\usepackage{graphicx}
\usepackage{color}

\newcommand{\comment}[1]{}

\newcommand{\Ll}{L}

\newcommand{\sCH}{{\sf CH}}
\newcommand{\CH}{\mbox{\sf CH}}
\newcommand{\Th}{T}
\newcommand{\KPU}{\mbox{\sf KPU}}
\newcommand{\KP}{\mbox{\sf KP}}
\newcommand{\ZF}{\mbox{\sf ZF}}
\newcommand{\Vv}{\mbox{\bf V}}
\newcommand{\Mm}{M}
\newcommand{\Fm}{\mbox{\it Fm}}
\newcommand{\Pow}{\mbox{$\mathcal P$}}

\newtheorem{thm}{Theorem}[section]
\newtheorem{defi}[thm]{Definition}

\newtheorem{prop}[thm]{Proposition}

\newtheorem{fact}[thm]{Fact}
\newtheorem{conj}[thm]{Conjecture}

\begin{document}

\title{Changing a Semantics: Opportunism or Courage?}
\author{Hajnal Andr\'{e}ka, Johan van Benthem, Nick Bezhanishvili, Istv\'an N\'{e}meti}
\thanks{This paper is dedicated to Leon Henkin, a deep logician who has changed the way we all work, while also being an always open, modest, and encouraging colleague and friend.}
\thanks{The authors wish to thank several colleagues for their helpful comments, including Solomon Feferman, Maria Manzano, Vaughan Pratt, Gabriel Sandu, and, especially, Jouko V\"{a}\"{a}n\"{a}nen.}
\thanks{The first and fourth authors would like to acknowledge the support of the Hungarian National Grant for Basic Research No T81188. The third author would like to acknowledge
the support of the Rustaveli Foundation of Georgia grant FR/489/5-105/11.}

\date{}

\maketitle


\begin{abstract} The generalized models for higher-order logics introduced by Leon Henkin, and their multiple offspring over the years, have become a standard tool in many areas of logic. Even so, discussion has persisted about their technical status, and perhaps even their conceptual legitimacy. This paper gives a systematic view of generalized model techniques, discusses what they mean in mathematical and philosophical terms, and presents a few technical themes and results about their role in algebraic representation, calibrating provability, lowering complexity, understanding fixed-point logics, and achieving set-theoretic absoluteness. We also show how thinking about Henkin{'}s approach to semantics of logical systems in this generality can yield new results, dispelling the impression of adhocness.
\end{abstract}

\section {General models for second-order logic}
In Henkin \cite{Henkin50} general models were introduced for second-order logic and type theory that restrict the ranges of available predicates to a designated family in the model. This enlarged model class for higher-order logics supported perspicuous completeness proofs for natural axiom systems on the pattern of the famous method in Henkin \cite{Henkin49} for proving completeness for first-order logic. To fix the historical background for this paper, we review the basic notions here, loosely following the exposition in van Benthem and Doets \cite{vBDoets83}.

The language of second-order logic has the usual vocabulary and syntactic constructions of first-order logic, including quantifiers over individual objects, plus variables for $n$-ary predicates, and a formation rule for existential and universal second-order quantifiers $\exists X\varphi, \forall X\varphi$. {\it Standard models} $M = (D, I)$ for this language have a domain of objects $D$, and an interpretation map $I$ for constant predicate symbols. \footnote{For simplicity only, we mostly omit individual constants and function terms in what follows.} Next, {\it assignments} $s$ send individual variables to objects in $D$, and predicate variables to real predicates over $D$, viewed as sets of tuples of objects. Then we have the following standard truth condition:

\medskip
\vspace{0.2cm}

    $M , s \models \exists X\varphi  \; \text{ iff  \; there is some predicate $P$ of suitable arity with} \; M , s[X:=P] \models \varphi $

\medskip
\vspace{0.2cm}

Now we take the control of available predicates in our own hands, instead of leaving their supply to set theory. A {\it general model} is a tuple $(M, V)$ where $M$ is a standard model and $V$ is some non-empty family of predicates on the domain $D$ of $M$. There may be some further constraints on what needs to be in the family, but we will look into this below.

Now interpretation on standard models is generalized as follows:

\medskip
\vspace{0.2cm}

$M , s \models \exists X\varphi \; \text{ iff \;  for some predicate $P$} \; in \; V \;  \text{of suitable arity:} \; M , s[X:=P] \models \varphi $

\medskip
\vspace{0.2cm}

Henkin \cite{Henkin50} proved the following seminal result:

\begin{thm}
Second-order logic is recursively axiomatizable over general models.
\end{thm}

This can be shown by adapting a Henkin-style completeness proof for first-order logic.

\medskip

However, another insightful road for obtaining  Thm.\ 1.1 is via translation into {\it many-sorted first-order logic}.  \footnote{In what follows in this section, merely for convenience, we shall deal with {\it monadic second-order logic} only, where second-order quantifiers run only over unary predicates, or sets.} Consider a first-order language with two sorts: `objects' and `predicates', each with their own set of variables. In addition, the domains are connected by a special predicate $E xp$ saying that object $x$ belongs to predicate $p$. Now we can translate the language of second-order logic into this two-sorted first-order language via a straightforward map $\tau$. The clauses for the logical operators are obvious, with $\exists x$ going into the object domain, and $\exists X$ into the predicate domain. At the atomic level, the translation is as follows:

$$  Xy  \; \text{   goes to} \;  E yX $$

\smallskip
\vspace{0.2cm}

 In addition, we state one principle that makes predicates in the two-sorted first-order language behave like real set-theoretic predicates in general models:

$$  Extensionality \; (EXT) \quad \quad \forall pq: \; (p=q \leftrightarrow \forall x(E xp \leftrightarrow E xq)) $$

 \smallskip

\begin{prop}
For all second-order formulas $\varphi$, $\varphi$ is valid on all general models iff
    $\tau(\varphi)$ follows from $EXT$ on all models for the two-sorted first-order language.
\end{prop}

The proof is obvious except for one observation. Interpretations for `predicate variables' $p$ in two-sorted first-order models $M$ can be arbitrary objects in the domain. But Extensionality allows us to identify such objects $p$ one-to-one with the sets $\{ x
 \mid E xp \; \text{in} \; M \}$.

 \medskip

As a result of this embedding, the validities of second-order logic over general models are axiomatizable, and by inspecting details, one can  extract an actual axiomatization close to that of two-sorted first-order logic. Moreover, this correspondence can be modulated. We can  drop Extensionality if we want even more generalized models, where we  give up the last remnant of set theory: namely, that predicates are identified with their set-theoretic extensions. This generalization  fits well with intensional views of predicates as concepts or procedures, that occur occasionally in logic and more often in philosophy.

 But more common in the logical literature is the opposite direction, an even stronger form of set-theoretic influence, where we turn further second-order validities into constraints on general models. A major example is that one often requires the following schema to hold in its two-sorted first-order version:

 \medskip

    {\it Comprehension}
    \begin{center}  $\exists X \forall y (Xy \leftrightarrow \varphi)$ \hspace{2mm}  for arbitrary second-order formulas $\varphi$ with $X$ not free in $\varphi$. \end{center}

 \medskip

What this achieves in the above completeness proofs is enforcing the validity of what many people find a very natural version of the logical law of existential instantiation: \footnote{Comprehension does make the strong philosophical assumption that logical constructions out of existing predicates deliver available predicates, something that might be debated.}

$$
\varphi(\psi) \to \exists X \varphi(X)
$$

\medskip

One can add many other constraints in this way without affecting completeness. One example is the object dual to Extensionality, $\forall xy(x= y \leftrightarrow \forall p (Exp \leftrightarrow Eyp))$
 that expresses a sort of `Individuality' for objects in terms of available predicates.

How far are these first-order versions of second-order logic removed from the logic over standard models? The distance can be measured by means of one standard formula --
as established by Montague \cite{Montague65}.

\begin{fact}
A second-order formula is valid on standard models iff its first-order translation
$\tau \varphi$ follows from EXT plus the second-order formula $\forall X \exists p \forall y( Xy \leftrightarrow E yp)$.

\end{fact}

The newly added axiom tells us that all set-theoretic predicates lie represented in the first-order model. Amongst other things, this observation shows that, taking prenex forms, the full complexity of second-order logic resides in its existential validities.

\smallskip

The preceding results turn out to extend to both classical and intensional type theories (Gallin \cite{Galin75}) and other systems. For a modern treatment of the current range of generalized models for higher-order logics and related systems, we refer to Manzano \cite{Manz96, Manz14}.  \footnote{There is much more to second-order logic than the general perspective given here. In particular,  deep results show that well-chosen fragments of second order logic, over well-chosen special classes of standard models, can have much lower complexity, and  have surprising combinatorial content, for instance, in terms of automata theory: cf. Gr\"adel, Thomas, Wilke, eds. \cite{GTW02}.}

\section{Clearing the ground: some objections and rebuttals}
Before going to concrete technical issues, we set the scene by listing some perspectives on generalized models that can be heard occasionally in the logical community, if not read in the literature. \footnote{In what follows, for brevity, we write `general models' instead of `generalized models'.} We do this mainly to clear away some rhetoric, before going to the real issues. For instance,  the very term `standard model' is already somewhat rhetorical, as it prejudges the issue of whether other models might also be natural, by depriving them of a neutral name, and making their pursuit `non-standard' or `deviant'. \footnote{Sometimes, a little dose of linguistics suffices to dispel this rhetoric, e.g., by calling `non-classical' logics `modern' logics, or `non-intended' models `serendipitous' models.}

\medskip

General models are often considered an ad-hoc device with little genuine content. This objection can be elaborated in several ways. One may hold that the natural  semantics for second-order logic consists only of those models that have the full power set of the individual domain for their unary predicates, and likewise for higher arities. General models are then a proof-generated device lacking independent motivation and yielding no new insights about second-order logic. And there can be more specific objections as well. For instance, imposing comprehension as a constraint on general models smacks of circularity. In order to define  the available predicates, we need to know what they are already, in order to understand the range of the second-order quantifiers in the comprehension principles. But objections can also run in another direction. The general model-based completeness proof for second-order logic tells us nothing new beyond what we already knew for first-order logi!
 c. In particular, one learns nothing that is specific to second-order validity.
\footnote{Contrast this with the sense of achievement, based on innovative proof techniques, when axiomatizing a piece of second order logic on a class of intended models, such as the monadic second-order logics of trees with successor relations as the `initial segment' relation in Rabin \cite{Rabin62}.}

Now there is something to these objections, but the question is how much. At least, one will always do well to also hear the case for the opposition. Here are a few preliminary considerations in favor of general models, and the proofs involving them.

\medskip

First of all, it is just a fact that in many natural settings, we do not want the full range of all set-theoretically available sets or predicates, as it would trivialize what we want to say. For instance, in natural language, when Russell wrote that `Napoleon had all the properties of a great general', he meant all relevant properties, not some trivial one like being a great general.
Likewise, when we grasp Leibniz' Principle that objects with the same properties are identical, we do not mean the triviality that object $x$ has the property of being identical to object $y$,
but rather think of significant properties. More generally, while it is true that the intuitive notion of property or predicate comes with a set, its extension or range of application, it would be unwarranted set-theoretic imperialism to convert this into the statement that every set is a property or predicate.\footnote{One example that shows this is the discussion in philosophy of propositions viewed as sets of worlds. Clearly, propositions correspond to sets of worlds where they are true, but it does not follow at all that every set of worlds must be the extension of some proposition. We return to this theme at the end of Section 4.}

The same consideration applies to mathematics: even there, set theory is not the norm. In geometry, important shapes correspond to sets of points, but definitely not every set of points is a natural geometrical shape. Hilbert 1899 mainly axiomatizes points, lines, and planes, instead of points and sets of points, and that for excellent mathematical reasons. Likewise, topology looks at open sets only, many theories of spatial objects use only convex sets, and one can mention many similar restrictions. Also, at a more abstract level, category theory does not say that all function spaces between objects must be full, the definition of a category allows us control over the available morphisms between objects, and that is precisely the reason for the elegance and sweep of the category-theoretic framework.

A third line of defense might be a counter-attack, turning the table on some of the earlier objections. Intended models for second-order logic provide us with a magical source of predicates that we tend to accept without having enquired into how they got there, and whether they should be there. Moreover, these entities also come with a magical notion of truth and validity, that seems to be there without us having to do any honest work in terms of analyzing proof principles that would constrain the sort of entity that we want. From this perspective, it is rather the general models that force us to be explicit about what we want and why. This point can be sharpened up. It is often thought that general models decrease the interest of proofs establishing their properties, but as we shall see later, this is a mistaken impression. There are also natural settings where general models make logical proof analysis and technical results more, rather than less sophisticated.

\medskip

This introductory discussion is general and inconclusive. Moreover, in the process of getting some ideological issues out of the way, we may even have introduced new ones.  For instance, despite appearances, our aim with this paper is not criticizing the use of set-theoretic notions per se. But there is a distinction to be made. The set-theoretic language is a lingua franca for much of logic and mathematics, and it facilitates formulation and communication in a way similar to the role of academic English. What one should be wary of, however, is the often implicit suggestion of a further commitment to very specific claims of set theory qua mathematical theory in our logical modelling. \footnote{For this unwarranted  commitment, consider the much more radical claim that using academic English carries a commitment to `British' or `American values'.} But where the precise border line between set-theoretic language and `set theory' is located seems a delicate matter, and one wher!
 e studying general models may in fact provide a better perspective.

We will now look at more systematic motivations for the use of general models in logic, Henkin-style and eventually also beyond, and assess their utility in more precise terms.

\section{Logical perspectives on controlling predicates}
For a start, we continue with the two themes introduced in the preceding section: controlling the available predicates, and the fine-structure of reasoning.

\subsection*{Proof theory} Logical analysis of an area of reasoning, even of mathematical proofs, seldom uses all set-theoretically available predicates -- but only ones that are definable in some way. And such proofs are crucial: despite the claimed intuitive nature of standard models, we hardly ever model-check assertions in them, since these models are too complex and mysterious for that. Instead, we prove that certain second-order formulas hold, say, in the natural numbers, by means of a mathematical argument appealing to accepted general principles. And these proofs will usually employ only very specific predicates, often ones that are definable in some sense, witness a wide literature on formalizing mathematical theories from Bishop \cite{Bish67} to Reverse Mathematics \cite{Sim10}.\footnote{This is more delicate with non-constructive principles like the Axiom of Choice, which we forego here.}  A survey of predicative proof theories with, for instance, highly restricte!
 d comprehension axioms, would go far beyond the confines of this paper, whose main slant is semantic, but we submit that such systems embody a reasoning practice that goes well with Henkin-style models.

\subsection*{Definable predicates} Restrictions to definable predicates are natural even in Henkin's first-order completeness proof itself that started our considerations. The Hilbert-Bernays Completeness Theorem says that, in order to obtain counter-examples to non-valid first-order formulas, models on the natural numbers suffice where the interpreting predicates are $\Delta_{0}^{2}$. Even so, not every restriction of predicates to some family of definable sets will be a good choice of a general model class for completeness. For instance, Mostowski \cite{Mostowski55}
showed that first-order logic is not complete when we restrict the interpreting predicates in models to be recursive, or even recursively enumerable. Vaught \cite{Vaught60} considered the related problem of complexity for the sentences true in all constructive models. \footnote{
The related result that Peano arithmetic has no non-standard recursive models is often called `Tennenbaum's Theorem', see, e.g., Kaye \cite{Kaye11}.}
Likewise, if we constrain general models too much, say, to only contain predicates that are first-order definable in the underlying first-order model, we may get harmful  complexity. In particular, Lindstr\"om \cite{Lind73} proved that, if second-order variables range over first-order definable predicates only, the system is non-axiomatizable.
Thus, while a restriction to definable predicates may be a motivation for using general models, there is no general guarantee that this move lowers the complexity of axiomatizability
for the logic. Having too many predicates gives us second-order logic, having too few can also lead to high complexity. We should be in between, and where that lies precisely may differ from case to case.


\subsection*{New structure: dualities}

 Here is another mathematical consideration in favor of controlling predicates. Doing so drastically may reveal important structure that we would not see otherwise. For instance, consider the frequently rediscovered weak theory of objects and types in `Chu Spaces' (Barwise \& Seligman \cite{BS93}, Pratt \cite{Pratt99}, Ganter \& Wille  \cite{GW99}). This theory treats objects and types on a par, resulting in a very appealing duality between their behavior\footnote{E.g., in Chu spaces, Extensionality for predicates is exactly Leibniz' Principle for objects.},   constrained by a notion of structural equivalence capturing the basic categorial notion of adjointness, which allows for model-theoretic analysis (van Benthem \cite{vB00}). All this elegant structure remains invisible to us in standard models for second-order logic. Similar points hold in terms of algebraic logic, the topic of  Section 4 below.

 \subsection*{New theorems: dependence and games}
It may be thought from the above that working with general models will decrease proof strength, so that, if anything, we lose theorems in this way. But the case of dualities was an instance where we regained new theorems at a higher level that just do not hold for the original version on standard models. Here is another example of this phenomenon, where general models  increase mathematical content.

Consider so-called IF logic of branching quantifiers, introduced in Henkin \cite{H59}, and taken further by Hintikka as a general study of independence in logic (cf.\ Hintikka  \& Sandu \cite{HS95}).
\footnote{See V\"a\"an\"anen \cite{Vaa07} for a parallel development of a rich system of `dependence logic'.}
Enderton \cite{Enderton70} and Walkoe \cite{Walkoe70} showed that branching quantifiers are the existential functional fragment of second-order logic, having a very high complexity of validity as we have seen in Section 1. In line with this, there has been little proof theory of IF logic, leaving the precise nature of reasoning with independent quantification a bit of a mystery. However, one can also start from a natural deduction analysis of dependent and independent quantifiers (Lopez Escobar \cite{Lopez1991}) and get an insightful proof system. The natural complete semantics for this proof system turns out to be general models whose available functions satisfy some simple closure conditions. But there is more: these models also represent a crucial move of independent interest.

As is well-known, like first-order logic, IF logic has a semantics
in terms of evaluation games where truth amounts to existence of a
winning strategy for the `verifier' in a game with imperfect
information of players about moves by others (Mann, Sandu \&
Sevenster \cite{MSS11}). This existential quantifier over strategies
ranges over all set-theoretically available strategies, viewed as
arrays of Skolem functions. Now the general models correspond to a
restriction on {\it available strategies}, a very intuitive move in
modeling games played by real agents, and this can drastically
change the structure of the game.   In particular, what can happen
is that a game that is `determined' in the sense that one of the two
players has a winning strategy, now becomes non-determined. But then
a very natural mathematical language extension suggests itself.
Non-determined IF games still have {\it probabilistic solutions} in
mixed strategies by the basic theorems of von Neumann and Nash.!
  And the equilibrium theory of these new games is much richer than that of the original games. \footnote{Hintikka himself suggested a restriction
to `definable strategies' for yet different reasons, but such a restriction might be a case of the above over-simplicity inducing high complexity after all.}

\subsection*{Conclusion} This concludes our first foray into concrete logical aspects of general models, showing that they raise delicate issues
of calibrating proof strength, levels of definability for predicates, and new laws possibly involving attractive language extensions. What we see here is that, with general models, weakness is at the same time wealth. We now turn to a number of major technical perspectives that we will consider in more detail.

\section{General models and algebraic representation}

Algebraic semantics is one of the oldest ways of modeling logic,
going back to the seminal work of Boole, de Morgan, and Schroeder.
In this section, we investigate the relation between algebraic
models and general models, largely using one case study: the
algebraic semantics of modal logic. We identify a number of general
issues, revolving around representation and completeness, and toward
the end, we show how the interplay of algebraic semantics and
general models is  very much alive today. Thus we find a motivation
for general models that is very different from the  proof-theoretic
concerns that were central in earlier sections. \footnote{ In
another context, namely in finite-variable fragments of first-order
logic, Henkin \cite{H67} introduced homomorphic images of the
Lindenbaum-Tarski algebras as generalized models, and he explained
why and how they can be considered as models. In \cite{H73}, Henkin
calls these generalized models  `algebraic models'. He proves a
completeness theorem and uses these algebraic models to show
unprovability of a given formula. So, the present section shows that
the ideas originated from Henkin  \cite{H67, H73} are alive and
flourishing.}


\subsection*{Algebraic completeness of modal logic}
The algebraic semantics of modal logic is given by modal algebras
$(A, \Diamond)$, with $A$ a Boolean algebra and $\Diamond$ a unary
operation on $A$ satisfying $\Diamond 0 = 0$ and $\Diamond (a\vee b)
=  \Diamond a \vee \Diamond b$. Every modal logic $L$ is complete
with respect to the `Lindenbaum-Tarski algebra' of all formulas
quotiented by  $L$. \footnote{That is, two formulas $\varphi$ and
$\psi$ are equivalent if $L\vdash \varphi \leftrightarrow \psi$.}
The Lindenbaum-Tarski algebra of a logic $L$ validates all and only
the theorems of $L$, and hence every modal logic is complete with
respect to a natural corresponding class of modal algebras.
\footnote{More generally there exists a lattice anti-isomorphism
between the lattice of normal modal logics and the lattice of
equationally defined classes (varieties) of modal algebras: cf.
\cite{BRV01, CZ97, Kra99} for details.}

\subsection*{General models via algebraic representation and categorial duality}

Modal algebras are related to model-theoretic structures via natural representations. These model theoretic structures, often called `general frames', may be viewed as Henkin models for modal logic generalizing the original Kripke frames.
A \emph{Kripke frame} is a pair  $(X,R)$ of a set $X$  and a binary relation $R$ . Each Kripke frame gives rise to a modal algebra whose domain is the powerset $\mathcal{P}(X)$ of $X$ with the Boolean operations
$\cap$, $\cup$, $()^c$, plus an operation $\Diamond_R$ on $\mathcal{P}(X)$ defined
by setting $\Diamond_R(S) = \{x\in X: \exists y \in S$ such that $x R y\}$ for each $S\subseteq X$.

However, going in the opposite direction, modal algebras do not normally induce standard frames. However, here is a widely used representation method. Given a modal algebra $(A, \Diamond)$, take the space $X_A$ of all {\it ultrafilters} on $A$ and define a relation $R_A$ as follows; $x R_A y $ iff $a\in y$ implies $\Diamond a\in x$ for each object $a\in A$. Moreover, let $\mathfrak{F}(X_A) = \{\varphi(a): a\in A\}$, where $\varphi(a) = \{x\in X_A: a\in x\}$. The latter is a special family of `good subsets' of the frame $(X_A, R_A)$ satisfying a number of natural closure conditions: it forms a Boolean algebra and it is closed under the natural operation for the modality. This is an instance of the following general notion. A \emph{general frame} $(X,R, \mathfrak{F}(X))$ is a triple such that $(X,R)$ is  Kripke frame and $\mathfrak{F}(X)$ is a subset of $\mathcal{P}(X)$ closed under $\Diamond_R$, that is, if  $S\in \mathfrak{F}(X)$, then $\Diamond_R(S)\in \mathfrak{F}(X)$. 
\footnote{In particular, $(\mathfrak{F}(X), \Diamond_R)$ is a
subalgebra of $(\mathcal{P}(X), \Diamond_R)$.} General frames are
natural models for a modal language since they provide denotations
for all formulas -- and of course, they are general models in
Henkin's sense. Moreover, they are the right choice. A basic theorem
by J\'onsson \& Tarski \cite{JT51} generalizes the Stone
representation for Boolean algebras to the modal realm:

\begin{thm} Each modal algebra is isomorphic to the modal algebra induced by the general frame of its ultrafilter representation.
\end{thm}

But still more can be said. The general frames produced by the J\'onsson-Tarski representation  satisfy a number of special conditions with a natural topological background. They are {\it descriptive} in the sense that different points are separated by available sets, non-accessibility is witnessed in the available predicates, and a natural compactness or `saturation' property holds for the available predicates. While not all Kripke frames are descriptive, the latter property does hold for many other model-theoretic structures.
\footnote{For more details on descriptive frames we refer to \cite{CZ97, BRV01, Kra99, SV88}.} 
This is interesting as an example where a natural class of general models satisfies third-order closure conditions different from the ones considered in the original Henkin models for second-order logic.



The correspondence between modal algebras and descriptive frames can be extended to a full correlation between two mathematical realms (see any of \cite{BRV01, CZ97, Kra99, SV88, Ven07} for details):

\begin{thm} There exists a categorical duality between (a) descriptive frames and definable `$p$-morphisms' (the natural semantic morphisms between descriptive frames that preserve modal theories), and (b) modal algebras and Boolean-modal homomorphisms.
\end{thm}

Thus, well-chosen classes of general models support rich category-theoretic dualities.

\subsection*{Completeness and incompleteness for standard relational models} While all this is true, standard frames do play a central role in a major area of modal logic, its completeness theory. The basic completeness results of modal logic say that, for well-known logics $L$, theoremhood coincides with validity in the standard frames satisfying the axioms of $L$. For instance, a formula is a theorem of modal $S4$ iff it is valid in all reflexive-transitive frames -- and there are many results of this kind for many modal logics in the literature. By contrast, via the above representation, algebraic completeness would only match theoremhood with validity in the class of general frames for the logic, which is restricted to valuations that take values in the admissible sets. What is going on here?

Consider the Lindenbaum-Tarski algebra for any modal logic $L$. In general the Kripke frame underlying the general frame representation of this algebra (often called the `canonical general frame' for the logic) need not be a model for the axioms of $L$ under all valuations. But in many special cases, it is. One general result of this form is the well-known Sahlqvist Theorem which states that  every modal logic that is axiomatised by `Sahlqvist formulas' -- having a special syntactic form whose details need not concern us here --  is  complete with respect to a first-order definable class of standard frames.

Results like this are often seen as improving the automatic completeness provided by general frames. But even so, they do not detract from the latter's importance. The  proof of the Sahlqvist Theorem depends essentially on showing that Sahlqvist formulas that hold on the canonical general frame for the logic also hold in its underlying standard frame. Thus, general models can be crucial as a vehicle for {\it transfer} to standard models.\footnote{There are also links here to {\it modal correspondence theory} where we study relational properties expressed by modal axioms on Kripke frames. Correspondence Theory assumes a different shape on general frames, though, for instance, many of its classical results still hold when we assume that the available propositions in general frames are closed under first-order definability. We do not pursue this model-theoretic theme in this paper, though it is definitely a case where introducing general models also poses some challenges to exi!
 sting theory on standard models of a logic.}

 \medskip

However, there are  limits to the preceding phenomena. Completeness
theorems on frame classes are not always obtainable in modal logic.
Famous `incompleteness theorems' from the 1970s onward have shown
that there are consistent modal logics $L$ and formulas  $\varphi$
such that $L\cup \{\varphi\}$ is consistent, but $\varphi$ is not
satisfiable on any Kripke frame for $L$. Concrete examples are
non-trivial, and several interesting ones may be found in
\cite{vB78, BRV01, CZ97}). Incidentally, the consistency is usually
proved by exhibiting a general frame where the logic is valid.
\footnote{For a concrete example, consider the tense logic of
\cite{Thom74} which has L\"ob's axiom for the past, making the
relation transitive and well-founded on Kripke frames, and the
McKinsey axiom for the future, which states that above every point,
there is a reflexive endpoint. Taken together, these requirements
are inconsistent, but they do hold on the general frame consistin!
 g of the natural numbers with only the finite and co-finite sets as available propositions.} In fact, frame incompleteness is the norm among modal logics, witness the remarkable classification results in \cite{Blo78}, while a modern exposition can be found in \cite{CZ97}. What all this says is that, despite appearances, it is completeness for general frames that underlies deduction in modal logic, while Kripke-frame completeness is a bonus in special cases.

\subsection*{Further representation methods}
Our discussion may have suggested that the algebraic defense of general models depends on one specific representation method for modal logic. But in fact, there is a much broader theory where similar points can be made. For instance, there are more general representations for distributive lattices and Heyting algebras (Priestley \cite{Pri70, Pri72}, Esakia \cite{Esa74}, Davey \& Priestley \cite{DP02}, \cite{CZ97}), using prime filters instead of ultrafilters to analyze intuitionistic logic and related systems. But new representation methods are appearing continuously. One recent example is the categorical duality between the categories of `de Vries algebras' (complete Boolean algebras with a special relation $\prec$  satisfying natural conditions)  and compact Hausdorff spaces: \cite{DeV62, Bez10}. Significantly, all our earlier modal themes return in this much broader mathematical setting.

Going beyond this, there is also a flourishing representation theory for substructural logics, cf. Dunn  \cite{Dunn91}, Gehrke, Dunn, Palmigiano \cite{DGP05},
Gehrke \cite{Gehrke06}, Marra \& Spada \cite{MS12},   Galatos \& Jipsen \cite{GJ14}. However, in this extended realm, it seems fair to say that major open problems remain such as finding `good' representation theorems for residuated lattices. But the very fact that this is considered a serious mathematical open problem illustrates the importance attached to finding Henkin-style model-theoretic structures matching the algebras.



The preceding themes are not exclusively mathematical concerns, they also play inside philosophy. For various conceptual reasons, Humberstone \cite{Humb81} and Holliday \cite{Hol14} have proposed replacing the possible worlds semantics for modal languages by structures of `possibilities' ordered by inclusion. To make this work, one adapts the usual truth definition to clauses for Boolean operations with an intuitionistic flavor, and interprets the modality by means of a suitable relation or function among possibilities. This new framework leads to an interesting  theory that even improves some features of classical modal logic, and it can be understood in terms of regular open sets in topological spaces: cf. \cite{vBBH15}. But of interest to us here is an issue of representation. Each standard possible worlds frame naturally induces a possibilities model defined by extending a relation on a frame into a function on its powerset  \cite{Hol14} (this is similar to a coalgebraic!
  perspective of modal logic \cite{Ven07}). But conversely, it is not true that each possibilities frame can be represented as coming from a Kripke frame in this manner. What is the proper comparison then between the two realms? The solution is again a move to general models: what should be compared are possibilities models and {\it generalized frames} based on a Boolean algebra of regular open sets of a topological space (see \cite{Hol14}  and \cite{vBBH15} for details of this construction).

\subsection*{Conclusion}  Algebraic models for logical systems are a natural match with general models in Henkin's sense, or at least in Henkin's style, where precise connections are provided by a large and growing body of representation theorems. Of course, this is a special take on the genesis and motivation for general models which comes with many interesting features of its own.  For instance, we have drawn attention to special conditions on the general models produced by representation methods, and to the ongoing challenges of representation theory, while we have also shown how in this realm, standard models become an interesting special case, which sometimes, but not always, suffices for completeness.

\section{   General models, lowering complexity, and core calculi}

Henkin's general models lower the complexity of second-order logic to that of the recursively axiomatizable first-order logic. \footnote{This is interesting historically, since initially, first-order logic was not the measure of all things. It was proposed only in Hilbert \& Ackermann \cite{HA28} as a well behaved fragment of higher-order logic,  for which G\"odel \cite{Goedel29} then provided the first completeness result.} Behind this move,
 we can discern a more general theme: that of lowering complexity of given logical
 systems.\footnote{In \cite{H67}, Henkin introduces generalized models for the finite-variable fragments $L_n$ of
 first-order logic, and he proves a completeness theorem with respect to these generalized models.
 This is also an example for lowering complexity, since it is known that $L_n$ cannot have a finitistic
 Hilbert-style strongly complete proof system for the standard models (see, e.g., \cite[Thm. 4.1.3]{HMT85}).
 Here is a further example of Henkin's positive thinking. Reacting to the negative result just quoted
 about $L_n$, he initiated jointly with Monk in \cite[Problem 1]{HM74} the so-called finitization problem which
 in turn fruited interesting and illuminating results e.g., by Venema, Stebletsova, Sain, Simon (for further
 details we refer to \cite{S00}).
}
  And with that theme, there is no reason to stop at first-order logic, whose notion of validity is recursively enumerable, but undecidable.
We will see (in Section 7) a similar lowering in complexity in connection with models of computation.
Continuing in what we see as Henkin{'}s spirit, though going far beyond  the basic semantic approach presented in Section 1,  could there be generalized models for first-order logic that make the system decidable?

This might sound like a purely technical interest, but one can add a deeper consideration. Van Benthem \cite{vB96} observes how many `standard' logical systems, despite their entrenched status, embody a choice of modeling a phenomenon that combines basic features with details of the particular formal framework chosen. For instance, first-order logic wants to be a core calculus of Boolean connectives and quantification. But in order to do this job, its semantics is couched in terms of set-theoretic notions. Are the latter a harmless medium (the `lingua franca' of our earlier discussion), or do these `wrappings' add extraneous complexity? In particular, is the undecidability of first-order logic a logical core feature of quantification, or a mathematical reflection of the set-theoretic modelling? In order to answer such questions, again, we need a more general modeling for strategic depth, and a talent for sniffing out unwarranted second-order or otherwise over-specific featur!
 es. In this section, we present one line, so-called `general assignment models' for first-order logic, as an instructive exercise in generalized model thinking.  \footnote{The ideas in this section
go back to N\'emeti \cite{Nemeti86}, Venema \cite{Ven92}, van Benthem \cite{vB96}, Andr\'eka, van Benthem \& N\'emeti \cite{AvBN98}, to which we refer for details. Here, N\'emeti's work, couched in terms of `cylindric relativized set algebras', was inspired by remarks of Henkin about obtaining positive results in algebraic logic, rather than the wave of counter-examples in earlier stages of the theory.
We note that investigating relativized cylindric algebras also originates from Leon Henkin, see e.g., Henkin \& Resek \cite{HenRes}.}

The core semantics of first-order logic works as follows. \footnote{As in Section 1, we will disregard individual constants and function symbols for convenience.}
Models $M = (D, I)$ consist of a non-empty set of objects $D$ plus an interpretation map $I$ assigning predicates over $D$ to predicate symbols in the language. Semantic interpretation involves models $M$, formulas $\varphi$, and assignments $s$ of objects in $D$ to individual variables of the language. A typical and basic example of the format is the truth condition for the existential quantifier:

$$
M , s \models \exists x \varphi\ \ \mbox{ iff \ \ there is some $d$ in $D$ with $M, s[x:= d] \models \varphi$}
$$

\medskip

While the underlying intuitive idea of existential quantification is uncontroversial, this technical clause involves two noteworthy assumptions in terms of supporting machinery. The auxiliary indices of evaluation are taken to be maps from variables to objects, and it is assumed that each such function is available as an assignment. Thus, not in its predicate structure, but in its supporting structure for interpretation, a first-order model comes with a hidden standard second-order object, viz. the full function space $D^{VAR}$. \footnote{Later on, we will see what it is more precisely that makes such spaces high-complexity inducing, viz. their geometrical `confluence' or `glueing' properties.}

\subsection*{A modal perspective}

 This choice is not entirely obvious intuitively, as the first-order language could be interpreted just as well on an abstract universe of states s allowing for binary transitions between indices, yielding a basic clause for modal logic:

$$
M , s \models \exists x\varphi\ \  \mbox{ iff\ \   there is some $s'$ with $sR_x sÕ$ such that $M , sÕ\models \varphi$}
$$

\medskip

We can think of the state space here as some independent
computational device that regulates the mechanism of interpretation
for our language, or its `access' to the model. Of course, standard
models are still around as a special case.

Treated in this modal way, first-order logic retains the essentials of its compositional interpretation, but its core laws do not reflect any set-theoretic specifics of assignment maps. Rather, they form the `minimal modal logic' K that already takes care of a large slice of basic reasoning with quantifiers, including its ubiquitous monotonicity and distribution laws. \footnote{The modal character
shows in that the quantifier $\exists x$ is really a labelled modality $\langle x\rangle$ now.}

 The minimal modal logic is decidable and perspicuous, while its meta-theory closely resembles that of first-order logic (Blackburn, de Rijke,  Venema \cite{BRV01}).

\subsection*{General assignment models}

 On top of this modal base system, the additional valid principles of first-order logic lie in layers of more specific assumptions on the modal state machinery approaching the full assignment spaces of standard models. Here is one basic level, using the set-theoretic view of assignments as functions from variables to individual objects, but without the existence assumption of having all functions around.

\begin{defi}
A general assignment model is a structure $(M, V)$ with $M$ a standard first-order model $(D, I)$, and $V$ a set of maps from individual variables to objects in $D$.
Interpretation of first-order formulas in these models is standard except for the following clause:

\smallskip

\begin{center}
$M , s \models  \exists x \varphi$\ \   iff\ \   there is some $d$ in $D$ such that $s[x:=d]\in V$  and $M , s[x:=d] \models \varphi$.
\end{center}
\end{defi}

By itself, this is just a technical move, but again there is an interesting interpretation. `Gaps' in the space of all assignments encode possible dependencies between variables. Suppose that we have an assignment $s$ and want to change the value of $x$. Perhaps the only variant that we have available for $s$ in the special family $V$ also changes its value for $y$, tying the two variables $x, y$ together. One can then think of the logic of general assignment models as first-order logic freed from its usual independence assumptions.

\begin{thm}
First-order logic over general assignment models is decidable.
\end{thm}

Things start changing with yet further axioms that express existence assumptions on the assignments that must be present. Their content can be brought out by standard modal frame correspondence techniques. We give one illustration, stated, for greater familiarity, in terms of the earlier abstract modal models for first-order logic.

\begin{fact}
A modal frame satisfies the axiom $\exists x \forall y Sxy \to \forall y \exists x Sxy$ for all valuations
    (i.e., all interpretations of atomic formulas) iff its transition relations $R_x$, $R_y$
    for the variables $x, y$ satisfy the following Church-Rosser confluence property
    for all states $s, t, u$: if $s R_x t$ and $s R_y u$, then there is a $v$ with $t R_y v$ and $u R_x v$.
\end{fact}

In general assignment models, the confluence property says that we can make the corresponding changes of values for variables concretely step by step: say, with
$t = s[x:=d]$,  $u = s[y:=e]$ and $v = s[x:=d][y:=e]$. This particular property is significant:

\begin{fact}
Adding confluence to general assignment logic makes the logic undecidable.
\end{fact}

The reason is that assignment models satisfying the stated property are rich, or rather, regular, enough to run standard proofs of the undecidability of first-order logic through its ability to encode undecidable tiling problems in a two-dimensional grid
\cite{Harel85}. Stated in another way, again we see that the undecidability has to do with mathematical, rather than logical content of the modelling: its ability to express regular geometrical patterns.

\subsection*{Richer languages}

 A striking new phenomenon that occurs with many forms of generalized semantics in logic is that it suggests richer languages to interpret, or at least, more sophisticated versions of the original logical language over standard models.

A concrete case for general assignment models are {\it polyadic quantifiers} $\exists x \varphi$ where $x$ is a finite tuple of variables. In standard first-order logic, this is just short-hand for iterated prefixes of quantifiers
$\exists x\exists y\dots \varphi$. On general assignment models, the natural interpretation for $\exists x \varphi$ is as existence of some available assignment in the model whose values on the variables in
$x$ can differ, a sort of `simultaneous re-assignment'. This is not equivalent to any stepwise iteration of single quantifiers. Nevertheless, the logic allows this:

\begin{thm}
First-order logic with added polyadic quantifiers is decidable
    over general assignment models.

\end{thm}

But language extensions can also increase complexity much more drastically. For instance, the general assignment models over a given first-order model $M$ form a natural family under extension of their assignment sets. But this is an interesting structure in its own right, and it makes sense to add a new extension modality interpreted as follows:

$$
M , V, s \models \Diamond \varphi\ \  \mbox{ iff\ \   there is some $V' \supseteq  V$ with $M, V', s \models \varphi$}
$$

\medskip

Now it is easy to see that standard first-order logic can be embedded in this richer language by interpreting a standard first-order quantifier $\exists x$ as a modal combination $\Diamond \exists x$.

This shows that the new logic with a modality across assignment sets is undecidable. We conjecture that this system is recursively enumerable, although we have not been able to find a straightforward argument to this effect. However this may be, our general point is that generalized models do not just lower complexity or provide axiomatizations for fixed logical formalisms. They may also change the whole design of languages and logics.

\subsection*{From general models to fragments}

 One proof of completeness and decidability for first-order logic over general assignment models mentioned above proceeds by translation into the decidable `Guarded Fragment' of first-order logic (GF, Andr\'eka, van Benthem \& N\'emeti
 \cite{AvBN98}). GF allows only quantifications of the following syntactic type: \footnote{We omit decidable extensions to `loosely guarded', `packed' and, recently, `unary negation' fragments.}

\begin{center}
$\exists y (G(x, y) \wedge \varphi(x,y))$, where $x, y$ are tuples of variables, $G$ is an atomic predicate, \\
and $x, y$ are the only free variables occurring in $\varphi (x, y)$.
\end{center}

\begin{thm}

 First-order logic over general assignment models can be translated faithfully
    into the Guarded Fragment of first-order logic over standard models.

\end{thm}

This is interesting, since a full language interpreted over generalized models now gets reduced to a syntactic fragment of that full language, interpreted over standard models. There is also a converse result (van Benthem \cite{vB05}) tightening the connection: the Guarded Fragment can be reduced to first-order logic over general assignment models.

Thus, again we see that general assignment models are a laboratory for rethinking what a generalized semantics means. Sometimes, a move to generalized models can also be viewed entirely differently, as one from a full logical language to a sublanguage, where new features of the generalized models show up as syntactic restrictions.

\subsection*{Conclusion} We have taken the spirit of general models one step further, and applied it to another `standard feature' of the semantics of logical systems, its use of assignments sending variables into the domain of objects. This method can be taken further than what we have shown here, since it applies to about any logical system. \footnote{For instance, while the base system of `dependence logic' in \cite{Vaa07} is undecidable and in fact of higher-order complexity, one can again find a standard power set in the background, the set of all sets of first-order assignments,  and once this is tamed, the core dependence logic will become decidable again.} We have found a significant border line between decidable core theories and  complexity arising from using set-theoretic objects such as full power sets. Moreover, as with algebraic representation, we found a number of further general issues that emerge in this special realm, such as possible redesign of logical la!
 nguages, and the interplay between generalizing semantics and passing to fragments.

%
%
%
%
%
%

\section{General models and absoluteness}

In this section, we discuss yet another general perspective on
general models, coming from the classical foundations of set theory
and logic itself -- and motivated, to a certain extent, by the
desire to draw a principled border line between the two. As with our
earlier topics, we cannot give an extensive technical overview, so
we merely present some basics along the way that are needed for the
points that we wish to make.\footnote{We thank Jouko V\"a\"an\"anen for sharing many useful corrections and insights concerning absoluteness, only some of which we have been able to do justice to in this brief paper.}

\subsection*{Absoluteness and non-absoluteness}

There is a sentence $\varphi_{\sCH}$ of second-order logic (SOL)
which expresses the Continuum Hypothesis (\CH) in the following
sense: for every model $\Mm$ with cardinality at least continuum we have
\begin{equation*}
\tag{*} \Mm\models\varphi_{\sCH} \quad\mbox{iff}\quad \CH\
\mbox{holds in our metatheory.}\footnote{Equivalently, \CH\ holds in the
`real world', cf.\ \cite[Sec. I.3, pp. 7-8]{Kunen80}.}
\end{equation*}

\noindent Hence $\varphi_{\sCH}$\footnote{More precisely, a suitable
variant of it taking into account the cardinality condition on $\Mm$. For the formula
$\varphi_{\sCH}$ see Ebbinghaus et al.\  \cite[pp.141-2]{EFT}.} is
logically valid in SOL iff $\CH$ holds in the meta-level set theory
`floating above our heads'. Thus, when contemplating the logical
validities of SOL, we have to rely on the set theory we use for
modeling/formalizing SOL, once more illustrating the
earlier-discussed phenomenon of `wrappings' versus `content'. Many
authors have agreed that the meaning of a formula $\varphi$ in a
model $\Mm$ should depend on $\varphi$ and $\Mm$ (and their parts)
and not on the entire set theoretic hyperstructure of the whole
universe containing $\varphi$ and $\Mm$, 
let alone on deciding highly complex and mysterious assertions like
the Continuum Hypothesis (cf. \cite{Kunen80}). And here is where
absoluteness kicks in. Meaningfully investigating these kinds of
questions, trying to make tangible and precise definitions and
statements concerning phenomena such as the above lead to the theory
of absolute logics.

A set-theoretic formula $\varphi$ is called {\it \ZF-absolute} if its truth
value does not change in passing from any model  \textbf{V} that
satisfies the axioms of standard Zermelo-Fraenkel set theory \ZF\ to
any transitive submodel $\textbf{V}'$ that satisfies those same
axioms. 
%
Next, Kripke-Platek set theory \KP\ can be considered as an austere
effective fragment of \ZF\ where we omit the Axiom of Infinity and we restrict
the Replacement and Collection schemes to bounded FOL-formulas.
\KPU\ is a natural extension of \KP\ where we also allow urelements. See
\cite[pp.10-12]{B75}. A set theoretic formula is said to be {\it
absolute} simpliciter if it is \KPU-absolute. A logic $\Ll=\langle \Fm,
\models\rangle$ is called {\it absolute} if the set-theoretic
formulas defining $\Fm$ and the ternary satisfiability relation
$\models$ are both absolute.


Here are a few illustrations for concreteness. Typical absolute
concepts are ``being an ordinal", ``being a finite ordinal", ``being
the $\omega$ (i.e., an ordinal with each element a successor ordinal,
but itself not a successor ordinal)", ``being the union of two sets
($x=y\cup z$)", etc. However, equally typically, the property ``$x$
is the powerset of $y$" is not $\ZF$-absolute, because whether a
collection of elements of $y$ is elevated to the rank of being a set
is independent from the $\in$-structure of $y$.

Intuitively, absoluteness of a logic means that truth or falsity of
the predicate $\Mm,s\models\varphi$ should depend only on the
$\in$-structures of $\varphi$, $\Mm$ and $s$, and not on the
`context' that $\Mm$, $\varphi$ and $s$ are living in (an example of
such a context would be the powerset of $\Mm$). Then our (*) above shows
that the logic SOL is not \ZF-absolute.

Barwise and Feferman \cite{B72} initiated the study of absolute
logics. The following theorem elucidating the special attention to
\KPU-absoluteness and making a crucial link to first-order logic is
due to Kenneth Manders \cite{M79}, see also Akkanen
\cite{Akkanen95}, Theorem 9 of Feferman \cite{F07} and Theorem 3.1.5
in \cite[pp. 620-622]{Vaa85}. Below, by $\Ll_{\omega\omega}$ we
understand many-sorted first-order logic (FOL) with equality. Also,
logic $\Ll'$ is said to be stronger than $\Ll$ if every class $C$ of
models defined by a sentence of $\Ll$ is also definable by a
sentence of $\Ll'$.

\begin{thm}\label{manders-thm}  $\Ll_{\omega\omega}$ is the strongest logic
among those absolute logics whose formulas are hereditarily finite
sets and whose structures $\Mm$ have domains consisting only of
urelements.
\end{thm}

Here we allow urelements in $\KPU$ so that we can have a good variety of
infinite models built from urelements without forcing the existence
of an infinite `pure set'.

\subsection*{General models via absolute versions}

If an important logic $\Ll$ turns out to be non-absolute, then it
seems useful to consider and study an absolute version\footnote{We
use the expression `absolute version' in an intuitive way, not in
the technical sense of \cite{B75}.} of $\Ll$ (besides $\Ll$ itself).
Indeed, Henkin's SOL is an absolute version of SOL. While SOL is not
$\ZF$-absolute, one can check that the satisfaction relation of
Henkin's general models for second-order logic is absolute (and it
even remains so if we add a finite number of our favorite axioms
valid in standard SOL). One can see this by using Manders' theorem
Thm. \ref{manders-thm} and by recourse to the results in Section 1,
in particular, the close connection between SOL on Henkin models and
many-sorted FOL. In the light of Manders' theorem, any absolute
version of SOL must have the form of a possibly many-sorted FOL
theory.




There is a general method for obtaining absolute versions of given
non-absolute logics $\Ll$. We may assume that the non-absolute\-ness
of $\Ll$ originates from the set-theoretical definition of the
ternary satisfiability relation of $\Ll$ not being absolute. Hence,
the set-theoretic formula defining the ternary $\Mm,s\models\varphi$
contains quantifiers ranging over some objects not in the transitive
closures of $\Mm,\varphi$, and $s$. Collect these `dangerous
objects' into an extra sort $S$. Consider the resulting extended
models $(\Mm,S,\in)$. Similarly to the example of SOL and its Henkin
models above, we arrive at a recursively axiomatizable many-sorted
FOL theory  $\Th$ whose `standard' models are the $(\Mm,S,\in)$'s. We
regard $\Th$ as an absolute version of the original logic $\Ll$. By
fine-tuning the axioms of $\Th$, we can fine-tune the absolute
version we want to work with. The models of $\Th$ which are not
standard, are the nonstandard models of the absolute version $\Th$
of $\Ll$. So, an absolute version of a logic always comes with a
notion of a generalized model for $\Ll$. Indeed, Henkin's models for
SOL can be seen as having arisen from devising an absolute version
of SOL, and later in Section 8 we will see other examples when the
notion of a generalized model for computation can be seen as
devising an absolute version of a non-absolute logic. Moreover, in
many cases a nonstandard model arising in this way can still
be seen as a standard one when viewed from another model of set
theory.

If we encounter a logic with high computational complexity, the
diagnosis may be non-absoluteness as the cause for high complexity.
This is the case with SOL and with First-Order Dynamic Logic with
standard semantics (see below). In all these cases, finding absolute
versions leads to natural notions of generalized models and to a
lower computational complexity of the absolute version. So, the
notions of the theory of absolute logics shed light on generalized
models as well as on lowering complexity.\footnote{While
absoluteness is meant to express independence from interpretations
of set theory, its very definition depends heavily on set theory --
but this may be considered a beneficial case of `catching thieves
with thieves'.}

At this high level of generality, the style of analysis in this
section applies very broadly. One of these applications concerns
logics of programs, or processes, cf.  \cite{ANSMFCS79, ANSTCS82}. Such logics
are usually \ZF-absolute, but not \KPU-absolute. They cannot have decidable proof systems,
even if we select acceptably small sublanguages in them, see
\cite[Thm.\ 1]{ANSMFCS79}. The reason for this is that standard
semantics assumes programs/processes to run in standard time, i.e.,
along the finite ordinals. We will return to this theme
in Sections 7, 8 below, where we consider generalized semantics for
recursion and computation in the setting of fixed-point logics.

\subsection*{Conclusion} Absoluteness as independence from set theory can be a
desideratum motivating forms of generalized semantics, and in fact, a
powerful methodology for their design.

\section{General models for recursion and computation}

Our final topic is an important phenomenon that has often presented
challenges to simple-minded semantic views of logical systems, the notion
of computation, and in particular, its characteristic nested structures of
recursion and inductive reasoning. These logics tend to be of very
high complexity, since treating computation explicitly often makes the natural
numbers definable by a formula of the logic, after which validity
starts encoding truth in second-order arithmetic.
\footnote{Likewise, so far, fixed-point logics of recursion, though
very natural at first glance, have successfully resisted
Lindstr\"om-style analysis in abstract model theory.}  

Intuitively,
logics of computation are about processes unfolding over time, and
hence it is crucial how we represent time: as a standard
set-theoretic object like the natural numbers, or as a more flexible
process parameter. We will start with this perspective, and its
generalized semantics. After that, we move to abstract fixed-point
logics of computation and generalized semantics that facilitate the
analysis of deductive properties and complexity-theoretic behavior.
In a final discussion, we show how two perspectives are related.

\subsection*{Nonstandard dynamic logics} Consider logics of programs in the tradition of propositional
dynamic logic $PDL$, based on program expressions with the regular
operations, where the key role is played by the Kleene star, i.e.,
{\it transitive closure}. While this propositional system is an
elegant decidable modal logic, its natural first-order version with
objects and predicates that applies to more realistic programming
languages is of highly intractable complexity. What causes this, and
how can this be remedied? We survey a few ideas from `nonstandard
dynamic logic' $NDL$ \footnote{The nonstandard models for this
dynamic logic were influenced by Henkin's generalized models for
SOL. More positively, one might also call these systems `logics of
general computation'.} \cite{ANSTCS82, CPZML84, GU91, GHR93, H81,
PTCS90}.

First-order dynamic logic cannot have a decidable proof system, even
if we select a small sublanguage, as is shown technically in
\cite{ANSMFCS79}.
The underlying reason is that standard semantics assumes processes to run in standard time, i.e., along the
natural numbers. If we make dependence on time explicit, however, we get a dynamic logic with an explicit time structure, or in other words, with a generalized semantics. In this nonstandard dynamic logic $NDL$, we are not tied to the natural numbers as the only possible time structure, but we can still make as many explicit requirements (as axioms) about  time as we wish. As a consequence of this greater flexibility in models, $NDL$ has a decidable proof system, see \cite{ANSTCS82}. \footnote{Here the proof
system is decidable, not the set of validities: Section 8 below has a question about this.}

More concretely, $NDL$ is a 3-sorted classical first-order logic
whose sorts are the time scale $T$, the data domain $D$, and a sort $I$
consisting of (not necessarily all) functions from $T$ into $D$. We
think of the elements of  $I$ as objects changing in time
(`intensions'). Typically, the content of a machine register changes
during an execution of a program, and this register is modeled as an
element of  $I$. In  standard models, $T$ is the set of natural numbers
and  $I$  is the set of all functions from $T$  to  $D$.  In  models of
$NDL$, these domains can be much more general, though we can impose
basic reasoning principles about computation such as induction
axioms for  $I$ (over specified kinds of formulas talking about time),
comprehension axioms for intensions, and axioms about the time sort
such as successor axioms, order axioms, as well as Presburger Arithmetic, or
even full Peano Arithmetic. Also, special axioms about the data
structure may be given according to the concrete application at
hand.

$NDL$  talks about programs, processes, and actions, and it
can express partial and total correctness plus many other
properties such as concurrency, nondeterminism, or
fairness. \footnote{$NDL$ has been used for characterizing the
`information content' of well-known program verification methods,
for comparing powers of program verification methods, as well as for
generating new ones \cite{ANSTCS82}.}


\subsection*{Digression: extensions to space-time} Similar methods have been applied to space-time
theories, such as the relativity theory of accelerated observers,
cf. \cite{MNSzFP06, SzPhD09}. There are good reasons for formalizing
parts of relativity theory in first-order logic FOL \cite{AMN06}.
And then another standard structure for time emerges. \footnote{This
kind of logical investigations of relativity theory were motivated
by Henkin's suggestion (to Andr\'eka and N\'emeti) for `leaving the
logic ghetto', i.e., for applying logic in other areas of science.}
Physical processes happen in real time, where the world-line of an
accelerated observer is an intensional entity whose spatial location
changes with time. It is customary to take the relevant time scale
here to be the real number line. But this can again be generalized: this time, 
to create a special relativity theory $AccRel$ of accelerated
observers where time structure becomes an explicit parameter subject
to suitable physically motivated axioms on the temporal order and
the world lines of test particles. While this is a `continuous'
rather than a `discrete' temporal setting, many of our earlier
points apply. For instance, instead of induction on the natural
numbers, one will now have axioms of Dedekind Continuity  expressing
that, if a physical property changes from holding to not holding
along a world-line, then there is a concrete point of time when the
change takes place. \footnote{Cf. \cite{BenthemLogicofTime} for
extensive model-theoretic discussion of such properties.}
Furthermore, one uses comprehension axioms to ensure existence of
enough world-lines for physical purposes, e.g., for being able to
select the inertial world-lines as those that locally maximize time.

We merely mention one striking outcome of this type of analysis with generalized models of time.  The famous Twin Paradox of relativity theory predicts that one of two twins who departs on a journey and undergoes acceleration will be younger upon her return than the other twin who stayed put. From a logical perspective, deriving this result turns out to involve having the real line as a time-structure. Using our more general models, it can be shown that the twin paradox cannot be derived by merely imposing temporal axioms, it is also essential to know how physical processes are related to time in terms of induction or continuity axioms about how properties of test
particles change along time.

\subsection*{Modal fixed-point logic on generalized models} Let us now move to a much more abstract view of induction and recursion over temporal or ordinal structures as embodied in modal fixed-point logics. In particular, consider the common running example of the modal {\it mu-calculus} which has become a powerful mathematical paradigm for the foundational study of sequential computation. \footnote{The contrast with the earlier $NDL$ approach may be stated in terms of moving from `operational' semantics to
`denotational' semantics of programs (cf. \cite{PTCS90}), but we
will not elaborate on this theme here.}

Consider the standard modal language enriched with a least fixed-point operator
$\mu x. \varphi$ for all formulas $\varphi$, where $x$ occurs under the scope of an even number of negations.
\footnote{For simplicity we disregard the greatest fixed-point operator, definable as $\nu x \varphi = \neg \mu x \neg \varphi$.} By a {\it general frame} for fixed point logic we mean a general frame $(X, R, \mathfrak{F}(X))$ as in Section 4 such that the family
$\mathfrak{F}(X)$ is closed under the fixed point operators. More precisely, here, inductively, each positive formula
$\varphi$ induces a monotone map $F_{\varphi}: \mathfrak{F}(X)\to \mathfrak{F}(X)$. \footnote{Here for $S\in \mathfrak{F}(X)$, we have that  $F_{\varphi}(S)$ is the value of $\varphi$ when
the variable  $x$ is assigned to $S$.}
Next we take the intersection of all pre-fixed points of $F_{\varphi}$, from $\mathfrak{F}(X)$, i.e., we consider
$\bigcap \{ S\in \mathfrak{F}(X)\ | \ F_{\varphi}(S)\subseteq S\}$. If this intersection belongs to
$\mathfrak{F}(X)$ for each $\varphi$,  then we call $(X, R, \mathfrak{F}(X))$ a \emph{general frame for modal fixed-point logic}.
This intersection is then the least fixed-point of $F_{\varphi}$ and it is in fact exactly the denotation of the fixed-point formula $\mu x. \varphi$.

\smallskip

This generalized semantics provides a new  way to interpret
fixed-point operators. Often, say, in spatial logic, we need to
restrict  attention to some practically realizable subsets of the
plane, and fixed-point operators need to be computed with respect to
these subsets only.  Note also that we lose nothing: if
$\mathfrak{F}(X) =  \mathcal{P}(X)$, then our generalized truth
condition coincides with the standard semantics of the modal
mu-calculus.

A major motivation for studying general-frame semantics for
fixed-point logic is that, via existing algebraic completeness and
representation methods, every axiomatic system in the language of
the modal mu-calculus is complete with respect to its general
frames.  Moreover, the powerful Sahlqvist completeness and
correspondence results from Section 4 extend from modal logic to
axiomatic systems in the modal mu-language for this semantics: cf.
\cite{BH12}. To appreciate this, we note that completeness results
for axiomatic systems of modal fixed-point logic with respect to the standard semantics
are very rare, and require highly complex machinery: see
\cite{Kozen83, Wal95}, and also \cite{SantVen10, tenCF10}. 

As a
further instance of the naturalness of these generalized models, we
note a delicate point of algebraic representation theory. Axiomatic systems of
the modal `conjugated mu-calculus' axiomatized by Sahlqvist formulas
are closed under the well-known {\it Dedekind-MacNeille completions}
in the above general-frame semantics (cf.\ \cite{BH12b}), whereas no
such result holds  for the standard semantics  (cf.\ \cite{Sant08}).


Now, following the general line in Section 4, we will quickly
overview the algebraic semantics for the modal mu-calculus. A modal
algebra $(A, \Diamond)$ is a \emph{modal mu-algebra} if for each
formula $\varphi$ positive in $x$,  the meet in $A$ of all the
pre-fixed points of $\varphi$ exists (see \cite{AKM95} and
\cite{BH12} for details). This meet will be exactly the denotation
of the fixed-point formula $\mu x. \varphi$. Similarly to modal
logic, the background for our earlier claims of completeness is that
axiomatic systems of the modal mu-calculus are  complete with
respect to modal mu-algebras obtained via the Lindenbaum-Tarski
construction: cf. \cite{AKM95, BH12}. A modification of the duality
between modal algebras and general frames  to modal mu-algebras
leads to the following notion. A descriptive general frame  $(X, R,
\mathfrak{F}(X))$ is called a \emph{descriptive mu-frame} if it
is a general mu-frame. For each descriptive mu-frame $(X, R,
\mathfrak{F}(X))$, the modal algebra  $(\mathfrak{F}(X_A),
\Diamond_{R_A})$ is a modal mu-algebra. Moreover, we have the
following converse: every modal mu-algebra   $(A, \Diamond)$ is
isomorphic to $(\mathfrak{F}(X_A), \Diamond_{R_A})$ for some
descriptive mu-frame $(X, R, \mathfrak{F}(X))$.


\smallskip

The difference between the standard and general-frame semantics for fixed-point operators shows in the following example.
Consider the frame $(\mathbb{N} \cup \{\infty\},R)$ drawn in Figure~1.
We assume $\mathfrak{F}(\mathbb{N}\cup \{\infty\}) =\{$all finite subsets of $\mathbb{N}$ and cofinite subsets containing the point $\infty\}$.
The standard semantics of the formula $\Diamond^\ast p$ is the set
of points that `see points in $p$ with respect to the reflexive and
transitive closure of the relation $R$'. Therefore, it is easy to
see that $\Diamond^\ast p$ is equal to the set $\mathbb{N}$. Indeed,
$\mathbb{N}$ is the least fixed-point of the map $S \mapsto \{0\}
\cup \Diamond_R (S)$.
However, if we are looking for a least fixed-point of this map in general-frame semantics, then we see that this will be the set $\mathbb{N}\cup \{\infty\}$.
We can say that the semantics of the formula $\Diamond^\ast p$ with respect to the general-frame semantics is the set of all points that  `see points in
$p$ with respect to the \emph{`admissible'  reflexive and transitive closure} of the relation $R$'.\footnote{A similar operation is used in \cite{Venema04} for characterising in dual terms subdirectly irreducible modal algebras.}

\begin{figure}
\begin{center}
\begin{tikzpicture}[->,thick]

\node (0) at ( -3,0) [circle,fill=black,inner sep=1pt,minimum size=8pt,label=above:$0$, label=below:$p$] {};
\node (1) at ( -1.8,0)  [circle,fill=black,inner sep=1pt,minimum size=8pt,label=above:$1$] {};
\node (2) at ( -0.6,0)  [circle,fill=black,inner sep=1pt,minimum size=8pt,label=above:$2$] {};
\node (infty) at ( 2.5,0) [circle,fill=black,inner sep=1pt,minimum size=8pt,label=above:$\infty$] {};
\node (d1) at ( 0,0) [circle,fill=black,inner sep=1pt,minimum size=2pt] {};
\node (d2) at ( 0.5,0) [circle,fill=black,inner sep=1pt,minimum size=2pt] {};
\node (d3) at ( 1,0) [circle,fill=black,inner sep=1pt,minimum size=2pt] {};
\draw [->,thick] (1) -- (0);
\draw [->,thick] (2) -- (1);
\path (infty) edge [loop right] (infty);
\end{tikzpicture}
\end{center}
\label{fig:ex}
\caption{Example}
\end{figure}

\medskip

While all this seems to be a  smooth extension from the basic modal case, there are also some deficiencies of this generalized semantics. We do not yet have many convincing concrete examples of general mu-frames, and especially of descriptive mu-frames. Some simple examples are general frames $(X, R, \mathfrak{F}(X))$, where $ \mathfrak{F}(X)$ is a complete lattice. For instance, the lattice of regular open sets $\mathcal{RO}(X)$ in any topological space is a complete lattice, so it allows interpretations of the least and greatest fixed-point operators. However, the meets and joins in such structures $\mathcal{RO}(X)$ are not set-theoretic intersection and union, and analyzing such concrete generalized models for their computational content seems to deserve more scrutiny than it has received so far.


\subsection*{Digression: a concrete fixed-point logic for recursive programs.} Modal fixed-point logics can be related to the earlier systems presented in this section. Our 3-sorted system $NDL$ essentially described operational semantics for flow-charts, whereas denotational semantics provides more structured meanings for programs than just input/output relations. In particular,  \cite{PTCS90} extends the Henkin-style explicit-time semantics of
$NDL$ to a Henkin-style denotational semantics for recursive programs (see also \cite[pp. 363--365]{GHR93}). The resulting system $NLP_{rec}$ is a 2-sorted classical first-order logic with one sort W for the set of states, and another
sort S for `admissible' binary relations over W, perhaps subject to Extensionality and Comprehension axioms. The semantics for recursive programs $p$ will now be a least admissible fixed-point in such a model.\footnote{E.g., take a  recursive program $p$ computing factorial: $p(x)   : =   \mbox{ if } x=0 \mbox{ then } 1
\mbox{ else } x\cdot p(x-1).$ The state-transition functional  $F$ associated with this program is $F( R )  :=  \{\langle
0,1\rangle\} \cup \{\langle x,x\cdot y\rangle : \langle
x-1,y\rangle\in R\}.$} (More details  will be discussed in Section 8.)

%

%

\subsection*{Comparing generalized dynamic logic and Henkin-style second order logic}
Both the logic $NDL$ and fixed-point logic over general models as
developed above are strongly connected to Henkin-style SOL (HSOL). We discuss one aspect of this connection here, though a full discussion is beyond the scope of this paper.

In one direction, $NDL$ may be viewed as a fragment of Henkin-style
SOL:

\begin{thm}
$NDL$ admits a natural translation into the language of the Henkin-style second-order logic.
\end{thm}

\noindent {\bf Proof idea.}
We would like to interpret $NDL$ (hence $T, D, I$) in HSOL. How can we
define the time-structure $(T, +, \ast)$ of $NDL$ from the concepts
available in HSOL?\footnote{Here we need  more
than monadic HSOL: e.g., binary or ternary relations for expressing addition $+$ on $T$. 
Therefore we need to  simulate binary and ternary relations in 
many-sorted HSOL. We simulate binary relations with a ternary
``elementhood" relation $E^2\subseteq D\times D\times V$ such that
$E^2xyX$ means that $x,y$ are ``objects", $X$ is a binary relation
and $Xxy$ holds (or, in other words, $\langle x,y\rangle\in X$).
Similarly, we use a four-place elementhood relation $E^3$ for
simulating ternary relations, etc. Thus our HSOL model in the
general case is of form $(D, V, E, E^2, E^3, ...)$, i.e. a structure of the sort considered in Section
1.}
%
Here comes the trick: In our HSOL structure $(D, V, E, E^2, E^3,
...)$ we can talk about second-order objects living in $V$. That is,
we can talk about sets, binary and ternary relations. So we
reconstruct $T$  as a pair of sets, $H\subseteq D$ and $f\subseteq
D\times D$ [i.e.,  $f$ is a binary relation] with $f: D\to D$  such
that $(H,f)$ behaves as the natural numbers with the successor
function. For example, we can postulate that there is an object $0\in H$
which generates the whole of $H$ via $f$. To see this, it is
sufficient  to note  that for all $H'\subseteq D$ if $0\in H'$ and
$H'$ is closed under $f$ then $H'\supseteq H$. From the axioms of
HSOL we can prove that this definition is correct, i.e., that
$(H,f)$ is unique up to isomorphism. The latter property is again
expressible in HSOL. We skip the rest of the proof which uses
similar ideas. \hfill $\Box$
\bigskip


A general experience reflected by the literature is that all program
properties expressed in HSOL so far were found to be expressible in
$NDL$. This motivates the following.

\begin{conj}
The program properties  expressible in the language of Henkin-style SOL can also be expressed
in $NDL$.
\end{conj}

We leave a similar detailed comparative discussion of expressive power in generalized frames for the mu-calculus to another occasion.


\subsection*{Conclusion} We have included logics for computation in this paper since these are often considered a challenge to general model techniques in logic. Indeed, we found two existing approaches to generalizing the semantics of logics with recursive definitions: one replacing standard mathematical structures by a wider class of models for their theory, leading to `generalized counting' on linear orders, the other modifying standard fixed-point semantics in terms of available predicates for approximations. While these generalized approaches to modeling computation seem prima facie different from the earlier ones presented in this paper, in the final analysis, they seem to fall under the general Henkin strategy after all.

\section{Relating different perspectives on generalized semantics}
We have discussed a number of general perspectives on Henkin's general models and further ways of generalizing semantics by defusing some of the invested set theory. Obvious further questions arise when thinking about connections between the approaches that we have discussed. We will only raise a few of these questions here, and our selection is not very systematic. Even so, the following discussion may help reinforce the interest of thinking about the phenomenon of generalized models in general terms.

\subsection*{Generalized fixed-point semantics and absoluteness} We have stated an objection of high complexity against prima facie attractive program logics such as first-order dynamic logic on its standard models. However, another charge would be lack of absoluteness. For a concrete example, consider the fixed-point logic of recursive programs in Pasztor~\cite{PTCS90}.

The meaning of a recursive program $p$ in a model $\Mm$ is defined as the least fixed-point of a functional $F$ associated with $p$ as a unary function over the set  $W$ of all possible states for $p$ in $\Mm$, i.e., over $\Pow(W\times W)$. As an illustration for the idea, consider computing the transitive closure of a relation $R$ as the
denotational semantics of a nondeterministic recursive program
\[ p(x)=y\quad\mbox{if either $R(x,y)$ or else\ }\exists
z[p(x)=z\land R(z,y)] .\] The functional $F$ associated with this
program is 
\[ F(X) = R\cup (X\circ R),\quad \mbox{where $\circ$ denotes relational
composition.}\] It can be seen that the least fixed-point of the
above $F$ is defined by
\[\begin{array}{ll}\Phi(x,y)\leftrightarrow &\exists\mbox{
sequence $s$ whose domain is a finite ordinal $n+1$ such
that}\\
&(x=s_0, y=s_{n}\land\forall i<n\, \langle s_i,s_{i+1}\rangle\in R).
\end{array}\]
This definition of transitive closure is $\ZF$-absolute. If
$\Vv, \Vv'$ are $\ZF$-models with $\Vv'$ a transitive submodel of $\Vv$,
then $\Vv$ and $\Vv'$ have the same finite ordinals since $\omega$
is $\ZF$-absolute and $\Vv'$ has to have an $\omega$ by the axiom
of infinity. However, the definition is not $\KPU$-absolute. Let $\Vv',\Vv\models\KPU$ where $\Vv'$ is a transitive submodel of $\Vv$
from which $\omega$ is omitted. This is possible because there is no Axiom of Infinity in \KPU\ as defined in \cite[pp.10-12]{B75}.
It may happen that $\Vv$ has more finite
ordinals than $\Vv'$. Hence the sequence $s$  in $\Phi$ might exist in $\Vv$ but not in $\Vv'$.
Therefore, for some concrete $x,y\in\Vv'$ we may have $\Phi(x,y)$ in
$\Vv$ but not in $\Vv'$.  This shows that $\Phi(x,y)$ is not a \KPU-absolute formula.


%


%

Now Pasztor~\cite{PTCS90} describes an absolute version of this semantics,
which is precisely the logic $NLP_{rec}$ for recursive programs discussed earlier in Section 7.
In particular, its choice of `admissible transition relations' is entirely analogous to
introducing a first-order sort for predicates as we did in Section 1. As a result,
$NLP_{rec}$ has a complete and decidable inference system. A similar analysis shows how our earlier logic $NDL$ can be viewed as the result of making the semantics of first-order dynamic logic absolute in a principled manner.

As a final issue in this connection, recall the earlier-mentioned feature that $NDL$ is a recursively axiomatizable but not a decidable theory. Can we make it truly decidable, in line with constructive set theories such as those studied in
\cite{COP01}, by merging its semantics with the general assignment models of Section 5? \footnote{A similar move toward creating decidable theories in physics might follow the lead of Section 7 in having only specially 
designated `space-time trajectories' available in models for physics.}

A more general issue arising here is how our earlier motivations for decidable semantics extend to the realm of computation. Our semantics of generalized counting still imports a substantial first-order theory of data structures for 
the computation to work on, and seems an interesting matter for even closer scrutiny. For instance, in line with Section 5, we could also define {\it guarded fragments} of second-order logic whose bounded quantification 
patterns match restrictions on available assignments of denotations to individual and predicate variables. \footnote{Here guard atoms will now be syntactically `third-order' (though this term loses much of its force in a many-sorted setting), 
while we must also allow names for designated objects and predicates. We leave technical details of such more radical decidability perspectives for another occasion.} 

\subsection*{Smoking out set theory in algebraic representation} In Section 4, we used the case of modal logic to show how general models are completely natural from the viewpoint of algebraic semantics for logical systems -- and even the only vehicle that will lead to insightful categorial dualities between algebras and model-theoretic structures. 

Even so, second-order logic can hide in unexpected places. For instance, note that the standard Stone Representation for Boolean Algebra, on which the modal representation is based, {\it itself employs} a non-absolute standard set-theoretic object, viz. the {\it set of all ultraflters} on a  given algebra. This shows in various side effects. One is that the models produced by this method tend to have overly large point sets as their domains. For instance, to represent a countable Boolean algebra, one needs only countable many points to make all necessary distinctions, whereas the Stone Space will have uncountably many
points.\footnote{ To see that a countable structure suffices, take a
countable elementary substructure of the Stone Space.}

But a further natural concern may be the heavy dependence of this
representation on set-theoretic assumptions like the existence of
ultrafilters extending given filters (the `Prime Ideal Theorem').
Indeed, van Benthem \cite{vB86} raised the question whether this
feature can be eliminated, and came up with what is essentially a
new `possibility semantics' for classical logic where representation takes
place in the universe of all consistent sets, by generalizing the
semantics of the basic logical operations to the latter
setting.\footnote{This theme will be resumed in the forthcoming
study \cite{vBBH15}.} 

This move suggests a further kind of generalized
semantics, where we do not just extend the model class, but where we
must also radically generalize our understanding of the meanings of
basic expressions in our original language. One way of thinking
about this is as turning the tables, using generalized semantics on
the very representation method that is supposed to motivate
generalized semantics. We leave this as an intriguing unexplored
avenue for now.

\subsection*{Richer languages for generalized models}

 We have seen that generalized assignment semantics supports richer first-order languages with natural new operators. This phenomenon is well-known from other areas of logic, such as linear logic and other resource- conscious logics, where classical logical operators `splitÕ into several variants, while new operators emerge as well. Again, an immediate test question is then whether the same happens in second-order logic, the realm where Henkin started his analysis. Do generalized predicate models support natural extensions of the usual second-order language?

We have not been able to find a good example that stays inside
FOL. However, one can add a new modality as we did for general
assignment models in Section 5, that says that the matrix formula is
true {\it in some extension of the current predicate range}.
Interesting things then become expressible, and the usual language
of second-order logic gets embedded, but we are not even sure what
sort of higher-order logic this system would really be.

\subsection*{Trade-offs between generalized semantics and language fragments}
We have seen a close analogy between first-order logic with general assignment models and the Guarded Fragment over standard models, as almost two sides of the same coin. Now well-chosen fragments are of great logical interest, and we noted their importance already in understanding how second-order logic works at deeper semantic and computational levels. What can we say about this in general? Can we re-analyze basic fragments of second-order logic in terms of generalized semantics for the whole language?  \footnote{As a related issue, could well-known key fragments of $SOL$ be captured through their insensitivity in passing between standard models and matching generalized models yet to be found?}


\subsection*{Coda: nonstandard truth conditions}
We conclude with one possibly more radical style of generalizing a semantics, by changing truth definitions, a theme that has come up occasionally in our earlier discussions, e.g., on constructivizing representation methods. 

As a concrete instance, consider the general assignment semantics of Section 5 as embodying our theme of `lowering complexity', taking FOL from undecidable to decidable by
removing extraneous mathematics from the core definitions of the logical system. Can we lower complexity even further, from decidable to `feasible' logics, whose satisfiability or validity problem would be decidable in polynomial time? No such proposal seems to exist, and it may be of interest to see where the difficulty lies. The complexity of decidability for many logics has to do with the exponential combinatorial explosion associated with {\it disjunctions}. Enlarging the model class does not seem to solve this, and what may be needed is rather a change in the way in which we interpret disjunction.   \footnote{One such proposal from the folklore has been to read disjunction as linear combination.}

Generalized semantics by tinkering with the interpretations of
logical operators can be very powerful. One interesting example of lowering complexity while improving model-theoretic behavior are
modal `bisimulation quantifiers' (Hollenberg \cite{Holle96}) that
step outside of current models in their interpretation. A
bisimulation quantifier is a modality $\langle p \rangle \varphi$
interpreted in pointed models $(M, s)$ as saying that there exists a
bisimulation, for the whole language except for the propositional
letter $p$, between $(M, s)$ and some model $(N, t)$ such that
$\varphi$  is true in $(N, t)$. One can think of this as a `tamed
form' of second-order quantification over the property $p$ -- and its effect is to
leave modal logic decidable and even preserves its striking
model-theoretic properties such as uniform
interpolation.\footnote{Added to propositional dynamic logic,
bisimulation quantifiers yield the modal mu-calculus.}   

A similar move would be possible for second-order logic, replacing
bisimulation by potential isomorphism, generalizing its semantics in
very different ways from the above. However, we will not pursue this
more radical form of meaning-changing generalized semantics here.

\section{Conclusion}
This paper started with an account of Henkin models for second-order
logic. To some people, this technique may look like a trick for lazy
people, but we have shown how this circle of ideas keeps returning
in the field for good reasons. We surveyed some major manifestations
of generalized semantics ranging from set theory to algebra and
computational logic, some close to Henkin's models, others more
remote -- and explored the many interesting technical and conceptual
issues and results that these bring to light.  \footnote{Even so,
our survey is by no means complete. E.g., Gabriel Sandu has
suggested that the `substitution account' for the meaning of
first-order quantifiers would be another natural candidate for our
style of analysis.} Our main suggestion has been that it is
well-worth thinking about the phenomenon of generalized semantics
more generally, both in terms of general insights and concrete
technical results when their main features are out side by side. We
feel that we have only scratched the surface of what might be a
deeper theoretical understanding of what is going on here.
\footnote{To mention just one more delicate point, it is argued in
\cite{HKV12}  that second-order logic with generalized semantics is
model-theoretically stronger than first-order logic because even
weak comprehension in the second order part brings about `strong
instability' in the model theoretic sense. In another vein,
\cite{HKV12} compares second-order logic with generalized semantics
(plus comprehension) to full second-order rather than to first-order
logic. What emerges then  is that, by forcing one obtains
model-theoretic results for full second-order logic, but with
generalized models one can eliminate forcing.}

While we have proposed taking a critical look at unwarranted `set-theoretic imperialism' and unquestioning acceptance of set-theoretic structures without a cost-benefit analysis, we have not advocated doing away with set theory altogether, since that would mean abandoning a lingua franca that has served the field of logic so well. Likewise, we have not proposed  abandoning second-order logic: it is attractive to look at behavior on standard models, and that attraction will not easily go away.
\footnote{We refer to \cite{Vaa12} for a state-of-the-art discussion of set theory side by side with second-order logic.} Finally, our tolerance for generalized semantics may seem almost all-embracing at times, but we certainly would not advocate trying to make general sense of every auxiliary move by logicians. Some practices in a field are just ad-hoc: and nothing is gained by pretending otherwise.

    Finally, however, we would like to emphasize the original motivation for writing this paper once more. It is always worth-while  to think about the achievements of Leon Henkin and what they mean in a broader perspective, and we have done so triggered by the various references to his work in this paper. In addition, we hope that the free-thinking and sometimes playful way in which we have done so reflects something of Henkin's  open-minded  personality that one would do well to emulate.

\bibliographystyle{plain}
\bibliography{Henkin}

\end{document}